\documentclass[12pt]{article}
\usepackage{amssymb,amsmath,theorem}

\allowdisplaybreaks

 \let\eps=\varepsilon
  \let\kappa=\varkappa
 
 \newcommand\pin{\kern.0833em}                    
 \newcommand\abs[1]{\left|#1\right|}              
 \newcommand\norm[1]{\left\|#1\right\|}           
 \newcommand\qed{\ifhmode\unskip\nobreak\fi\quad  
    \ifmmode\square\else\hbox{$\square$}\fi}      
 \newcommand\proofskip{\vspace{
        \theorempostskipamount}}                  

 \newtheorem{theorem}{Theorem}
 \newtheorem{lemma}[theorem]{Lemma}
 {\theorembodyfont{\normalfont}\newtheorem{rem}[theorem]{Remark}
    \newtheorem{defn}[theorem]{Definition}
    \newtheorem{I}[theorem]{\unskip}}
 \newcommand\proof[1]{\noindent\textit{Proof#1}}

\begin{document}

 \vspace*{10mm}

 \begin{center}
 {\large\bfseries
 \uppercase{Variational principles for $\boldsymbol{t}$-entropy,\\[2pt]
 the spectral potential of transfer operator,\\[2pt] and entropy statistic theorem\\[2pt]
 are equivalent\\}}

 \bigskip \medskip

 {\large\mdseries\itshape
 V.\,I.\ Bakhtin $($bakhtin@tut.by$)$$^\dag$, \
 A.\,V.\ Lebedev $($lebedev@bsu.by$)$$^\ddag$}

 \bigskip \medskip

 $^\dag$Belarusian State University / John Paul II Catholic University of Lublin, Poland;\\
 $^\ddag$\,Belarusian State University / University of Bialystok, Poland
 \end{center}

\vspace{-18pt}

\renewcommand{\abstractname}{}

\begin{abstract} \noindent
For any transfer operator we establish the equivalence of variational principles for $t$-entropy,
the spectral potential and entropy statistic theorem and give new proofs for all these statements.
\end{abstract}

\medbreak


\medbreak


\bigskip\medskip


For any transfer operator we will establish the equivalence of the next statements: `variational
principle for $t$-entropy' $\Leftrightarrow$ `variational principle for the spectral potential'
$\Leftrightarrow$ `entropy statistic theorem'.

\section{Variational principles for $\boldsymbol t$-entropy and\\ the spectral potential of transfer
operators}

Let us start with recalling the main objects and notions in question.

\begin{I} \label{0..1}
We will call $\mathcal C$ a \emph{base algebra} if it is a self-adjoint part of a certain commutative
\hbox{$C^*$-al}\-ge\-bra with an identity $\bf 1$. This means that there exists a commutative
$C^*$-algebra $\mathcal B$ with an identity $\bf 1$ such that
\begin{equation*}
 \mathcal C= \{\pin b\in \mathcal B\mid b^* = b\pin\}.
\end{equation*}

\medskip

As is known the Gelfand transform establishes an isomorphism between $\mathcal C$ and the algebra
$C(X)$ of continuous real-valued functions on a Hausdorff compact space $X$, which is the maximal
ideal space of the algebra $\mathcal C$. Throughout the article we identify $\mathcal C$ with
$C(X)$ mentioned above.
\end{I}

The next known result (see, for example, \cite{end-R}) establishes a correspondence between
endomorphisms of base algebras and dynamical systems.

\begin{theorem} \label{6..2}
If\/ $\delta\!: \mathcal C \to \mathcal C$ is an endomorphism of a base algebra\/ $\mathcal C$ then
there exists an open-closed subset\/ $Y\subset X$ and a continuous mapping\/ $\alpha\!:Y\to X$
$($both\/ $Y$ and\/ $\alpha$ are uniquely defined\/$)$ such that
\begin{equation*}
 \big[\delta f\big](x) = \chi_{Y}(x)f(\alpha (x)), \qquad f\in {\mathcal C}, \quad x\in X,
\end{equation*}
where\/ $\chi_{Y}$ is the index function of\/ $Y$. In particular, if\/ $\delta({\mathbf 1}) =
{\mathbf 1}$ then\/ $Y = X$ and
\begin{equation} \label{6,,1}
 \big[\delta f\big](x) = f(\alpha (x)).
\end{equation}
\end{theorem}

\begin{rem} \label{6..3}
It is clear that any endomorphism of a $C^*$-algebra $\mathcal B$ is completely defined by its
restriction onto the self-adjoint part $\mathcal C$ of $\mathcal B$ and on the other hand any
endomorphism of $\mathcal C$ extends uniquely up to an endomorphism of $\mathcal B$. Therefore the
correspondence between endomorphisms and dynamical systems presented in the theorem can be equally
described in terms of endomorphisms of~$\mathcal B$.
\end{rem}

\begin{I} \label{6..4}
In what follows the pair $(\mathcal C,\delta)$, where $\mathcal C$ is a base algebra and $\delta$
is its certain endomorphism such that $\delta(\mathbf 1) = \mathbf 1$, will be called a
$C^*$-\emph{dynamical system}, and the pair $(X, \alpha)$ described in Theorem \ref{6..2} will be
called the \emph{dynamical system corresponding to} $(\mathcal C,\delta)$. The algebra $\mathcal C$
will be also called the \emph{base algebra of the dynamical system} $(X,\alpha)$.

Throughout the paper notation $\mathcal C$, $\delta$, $X$, $\alpha$ will denote the objects
introduced above and we will use either of them (say $\delta$ or $\alpha$) for convenience reasons
(once $\alpha $ is chosen then $\delta$ is defined uniquely by (\ref{6,,1}) and vice versa).
\end{I}

\begin{defn} \label{6..5}
Let $({\mathcal C},\delta)$ be a $C^*$-dynamical system. A linear operator $A\!:\mathcal C\to \mathcal C$ is
called a \emph{transfer operator}, if it possesses the following two properties

\smallskip

a) $A$ is positive (it maps nonnegative elements of $\mathcal C$ into nonnegative ones);

\smallskip

b) it satisfies the \emph{homological identity}
\begin{equation*}
 A\bigl((\delta f)g\bigr) =f Ag \quad \text{for all}\ \ f,g\in \mathcal C.
\end{equation*}
If in addition this operator maps $\mathbf 1$ into $\mathbf 1$ we will call it a \emph{conditional
expectation operator.}
\end{defn}

\begin{rem}\label{6..6}
Any transfer operator $A\!:\mathcal C\to \mathcal C$ can be naturally extended up to a transfer
operator on ${\mathcal B} = {\mathcal C} + i\pin{\mathcal C}$ by means of the formula
\begin{equation*}
 A(f+ig)= Af + i Ag.
\end{equation*}
On the other hand, given any transfer operator on $\mathcal B$, its restriction to $\mathcal C$
(which is well defined in view of property a) of Definition \ref{6..5}) is also a transfer
operator. Therefore transfer operators can be equivalently introduced as by means of $C^*$-algebra
$\mathcal B$ so also by means of its self-adjoint part --- the base algebra $\mathcal C$. We prefer
to exploit the base algebra since in what follows we use the Legendre transform which is an
essentially real-valued object.
\end{rem}

\begin{I} \label{6..8}
Let $(\mathcal C,\delta)$ be a $\mathcal C^*$-dynamical system and $(X,\alpha)$ be the
corresponding dynamical system. We denote by $M(\mathcal C)$ the set of all positive normalized
linear functionals on~$\mathcal C$ (which take nonnegative values on nonnegative elements and are
equal to $1$ on the unit). Since we are identifying $\mathcal C$ and $ C(X)$, the Riesz theorem
implies that the set $M(\mathcal C)$ can be identified with the set of all regular Borel
probability measures on~$X$ and the identification is established by means of the formula
\begin{equation*}
 \mu[\varphi] = \int_X \varphi\, d\mu, \qquad \varphi \in\mathcal C=C(X),
\end{equation*}
where $\mu$ in the right-hand part is a measure on $X$ assigned to the functional $\mu\in
M(\mathcal C)$ in the left-hand part. That is why with a slight abuse of language we will call
elements of $M(\mathcal C)$ \emph{measures}.

A measure $\mu\in M(\mathcal C)$ is called \emph{$\delta$-in\-va\-ri\-ant} if for each $f\in
\mathcal C$ we have $\mu[f]=\mu[\delta f]$. The set of all $\delta$-in\-va\-ri\-ant measures from
$M(\mathcal C)$ will be denoted by $M_\delta(\mathcal C)$. Clearly, in terms of the dynamical
system $(X,\alpha)$ the condition $\mu[f] =\mu[\delta f]$ is equivalent to the condition $\mu[f]
=\mu[f\circ\alpha]$,\ \ $f\in C(X)$. Therefore $M_\delta(\mathcal C)$ can be identified with the
set of all $\alpha$-invariant regular Borel probability measures on X.
\end{I}

\begin{I}\label{6..9}
Let $A\!: \mathcal C\to \mathcal C$ be a fixed transfer operator for a $C^*$-dynamical
system~$({\mathcal C}, \delta)$. In what follows we consider the family of operators $A_\varphi\!:
\mathcal C\to \mathcal C$, where $\varphi \in {\mathcal C}$, defined by means of the formula
$A_\varphi f := A(e^\varphi f)$. Evidently, all the operators in this family are transfer operators
for $({\mathcal C}, \delta)$ as well.
\end{I}

Let $\lambda(\varphi)$ be the logarithm of the spectral radius of $A_\varphi$:
\begin{equation} \label{6,,6}
 \lambda(\varphi) = \lim_{n\to\infty}\frac{1}{n}\ln \big\|A_\varphi^n\mathbf 1\big\|,\qquad
 \varphi\in \mathcal C\, .
\end{equation}
The functional $\lambda(\varphi)$ will be called the \emph{spectral potential} of transfer operator
$A$.

By a \emph{partition of unity} in the algebra $\mathcal C$ we mean any finite set $D =
\{g_1,\dotsс,g_k\}$ consisting of nonnegative elements $g_i\in \mathcal C$ satisfying the identity
$g_1+ \dots +g_k = \mathbf 1$.

Definition of $t$-entropy $\tau(\mu) $ is given in the following way:
\begin{gather} \label{6,,5}
 \tau(\mu) := \inf_{n\in\mathbb N}\frac{\tau_n(\mu)}{n}\,,\qquad
 \tau_n(\mu) := \inf_D\tau_n(\mu,D),\\[6pt] \label{6,,7}
 \tau_n(\mu,D) := \sup_{m\in M(\mathcal C)}\sum_{g\in D}\mu[g]\ln\frac{m[A^ng]}{\mu[g]}\,,\qquad
 \mu\in M_\delta(\mathcal C).
\end{gather}
The infimum in \eqref{6,,5} is taken over all the partitions of unity $D$ in the algebra $\mathcal
C$.

If we have $\mu[g] = 0$ for a certain $g\in D$, then we set the corresponding summand
in~\eqref{6,,7} to be zero independently of the value~$m[A^ng]$. And if there exists an element
$g\in D$ such that $A^ng = 0$ and simultaneously $\mu[g]>0$, then we set $\tau(\mu) = -\infty$.

\begin{theorem}[variational principle for $\boldsymbol{t}$-entropy] \label{..1}
Let\/ $(\mathcal C,\delta)$ be a\/ $C^*$-dyna\-mic\-al system, $A\!:\mathcal C\to \mathcal C$ be a
certain transfer operator for\/ $(\mathcal C,\delta)$, and\/ $A_\varphi =A(e^\varphi\,\cdot\,)$ for
all\/ $\varphi\in\mathcal C$. Then the following equality takes place\/$:$
\begin{equation} \label{,,6}
 \tau(\mu) =\inf_{\varphi\in \mathcal C} \bigl(\lambda(\varphi) -\mu[\varphi]\bigr), \qquad
 \mu\in M_\delta(\mathcal C)
\end{equation}
\emph{(where $M_\delta(\mathcal C)$ is the set of all positive normalized $\delta$-invariant linear
functionals on\/ $\mathcal C$).}
\end{theorem}

\begin{theorem} \label{..2}
For each linear functional\/ $\mu$ on\/ $\mathcal C$ that does not belong to\/ $M_\delta(\mathcal
C)$ the following equality takes place\/$:$\hbox{\quad}
\begin{equation} \label{,,7}
 \inf_{\varphi\in \mathcal C} \bigl(\lambda(\varphi) -\mu[\varphi]\bigr) = -\infty.
\end{equation}
\end{theorem}

These two theorems show that it is natural to put
\begin{equation} \label{,,8}
 \tau(\mu) =-\infty \quad \text{for all}\ \ \mu\in \mathcal C\rule{0ex}{1.7ex}^*\setminus
 M_\delta(\mathcal C).
\end{equation}
Then formulae \eqref{,,6} and \eqref{,,7} are united into one:
\begin{equation} \label{,,9}
 \tau(\mu) =\inf_{\varphi\in \mathcal C} \bigl(\lambda(\varphi) -\mu[\varphi]\bigr), \qquad
 \mu\in \mathcal C\rule{0ex}{1.7ex}^*.
\end{equation}

The proof of Theorems~\ref{..1}~and~\ref{..2} will be implemented in a number of steps.

\begin{lemma} \label{..3}
For any\/ $\varphi \in \mathcal C$ and\/ $\mu\in M_\delta(\mathcal C)$ one has
\begin{equation} \label{,,11}
 \lambda(\varphi) \ge \mu[\varphi] +\tau(\mu).
\end{equation}
\end{lemma}

\smallskip

\proof. Let us show that for any $\varphi \in \mathcal C$, $\mu\in M_\delta(\mathcal C)$,
$n\in\mathbb N$ and $\eps >0$ there exists a partition of unity $D$ such that
\begin{equation} \label{,,10}
 \eps + \frac{\ln\|A_\varphi^n\|}{n} \ge \mu[\varphi] +\frac{\tau_n(\mu,D)}{n}\,.
\end{equation}
Once this is done then by arbitrariness of $\eps >0$ inequality \eqref{,,11} follows from
\eqref{,,10} by taking infimum with respect to $D,\,n$.

So it is enough to verify \eqref{,,10}.

Let us introduce the notation
\begin{equation*}
 \label{e-Sn} S_n\varphi :=\varphi+\delta\varphi+\,\cdots\,+\delta^{n-1}\varphi, \qquad
 \varphi\in \mathcal C.
\end{equation*}
Applying $n$ times the homological identity to the operator $A_\varphi^n =(Ae^{\varphi})^n$, we
obtain
\begin{equation} \label{2,,3}
 A_\varphi^nf =A(e^\varphi A(e^\varphi\dotsm A(e^\varphi f)...)) =A^n\bigl(e^{S_n\varphi}f\bigr).
\end{equation}

For arbitrary numbers $n\in\mathbb N$ and $\eps>0$ we choose a partition of unity~$D$ such that on
the support of each function $g\in D$ the oscillation of function $S_n\varphi$ does not
exceed~$\eps$. This~$D$ is in fact the desired partition.

Set
\begin{equation*}
S_n\varphi(g) := \sup\{\pin S_n\varphi(x)\mid g(x)\ne 0\pin\},
\end{equation*}

 \medskip\noindent
where we identify $g$ with the corresponding function in $C(X)$ (cf. \ref{0..1} and \ref{6..4}).

Equality \eqref{2,,3} and concavity of the logarithm function imply the following inequalities for
all functionals $m\in M(\mathcal C)$ and $\mu \in M_\delta(\mathcal C)$:
\begin{gather*}
 \eps+\ln \| A_\varphi^n \| \,=\, \eps+\ln \| A_\varphi^n\mathbf 1\| \,\ge\,
 \eps+\ln m[A_\varphi^n\mathbf 1]
 \,=\, \eps +\ln\sum_{g\in D}m\bigl[A^n\bigl(e^{S_n\varphi}g\bigr)\bigr]\\[3pt]
 \,\ge\, \ln\sum_{g\in D}e^{S_n\varphi(g)}m[A^ng]
 \,\ge\, \ln\sum_{\mu[g]\ne 0}\mu[g]\frac{e^{S_n\varphi(g)}m[A^ng]}{\mu[g]} \displaybreak[2]\\[3pt]
 \,\ge\, \sum_{\mu[g]\ne 0}\mu[g]\ln\frac{e^{S_n\varphi(g)}m[A^ng]}{\mu[g]}
 \,=\, \sum_{\mu[g]\ne 0}\mu[gS_n\varphi (g)] +\sum_{\mu[g]\ne 0}\mu[g]\ln\frac{m[A^ng]}{\mu[g]}\\[5pt]
 \,\ge\, \sum_{\mu[g]\ne 0}\mu[gS_n\varphi] + \sum_{\mu[g]\ne 0}\mu[g]\ln\frac{m[A^ng]}{\mu[g]}
 \,=\,\mu[S_n\varphi] + \sum_{g\in D}\mu[g]\ln\frac{m[A^ng]}{\mu[g]}.
\end{gather*}
Passing in these inequalities to the supremum over $m\in M (\mathcal C)$ and taking into
account~\eqref{6,,7}, one obtains the inequality
\begin{equation*}
 \eps+ \ln\|A_\varphi^n\|\ge \mu[S_n\varphi] +\tau_n(\mu,D) =n\mu[\varphi] +\tau_n(\mu,D),
\end{equation*}
which implies \eqref{,,10}. \qed

\proofskip

To finish the proof of Theorem~\ref{..1} we need two more lemmas.

In the next lemma the notation $\lambda(\varphi,A)$ has the same meaning as $\lambda(\varphi)$ and
$\lambda(n\varphi,A^n)$ denotes the logarithm of spectral radius of the operator
$A^n(e^{n\varphi}\,\cdot\,)$.
\begin{lemma} \label{..4}
The following inequality takes place\/$:$
\begin{equation} \label{,,12}
 n\lambda(\varphi,A) \le \lambda(n\varphi,A^n), \qquad n\in \mathbb N.
\end{equation}
\end{lemma}

\proof. Note that for any natural $k$ one has
\begin{equation} \label{e.Snk}
 \exp\{S_{nk}\varphi\} =\pin
 \exp\biggl\{\sum_{i=0}^{n-1}\sum_{j=0}^{k-1} \delta^{i+nj}(\varphi )\biggr\} \pin=\,
 \prod_{i=0}^{n-1} \exp\biggl\{\sum_{j=0}^{k-1} \delta^{i+nj}(\varphi )\biggr\}.
\end{equation}
Let $c =\|\varphi\|$. Bearing in mind observations \eqref{2,,3}, \eqref{e.Snk}, and exploiting
H\"older inequality in the form
\begin{equation*}
 \nu[\psi_1 \dotsm \psi_n] \,\le\, \prod_{i=1}^{n} \nu \big[\abs{\psi_i}^n\big]^{1/n},
\end{equation*}

 \smallskip\noindent
where $\nu[\,\cdot\,] = m\big[A^{n(k+1)}(\,\cdot\,)\big]$ and $m\in M(\mathcal C)$, one obtains
\begin{gather*}
 e^{-nc}\pin m\big[A_\varphi^{n(k+1)}\mathbf 1\big] \,=\,
 e^{-nc}\pin m\big[A^{n(k+1)}\big( e^{S_{n(k+1)}\varphi}\mathbf 1\big) \big] \,\le\,
 m\big[A^{n(k+1)}\big( e^{S_{nk}\varphi}\mathbf 1\big) \big] \\[6pt]
 \,=\, \nu\bigg[\pin\prod_{i=0}^{n-1} \exp\!\pin\bigg\{\sum_{j=0}^{k-1} \delta^{i+nj}\varphi\bigg\}\bigg]
 \,\le\, \prod_{i=0}^{n-1} \nu \bigg[\exp\!\pin\bigg\{\sum_{j=0}^{k-1}
 \delta^{i}\big(\delta^{nj}(n\varphi )\big)\bigg\} \bigg]^{1/n}\\[6pt]
 \,=\, \prod_{i=0}^{n-1} m\bigg[A^{n-i}A^{nk}A^i\bigg(\!\exp\!\pin\bigg\{\sum_{j=0}^{k-1}
 \delta^{i}\big(\delta^{nj}(n\varphi)\big)\bigg\}\pin \mathbf 1\bigg)\bigg]^{1/n}\\[6pt]
 \,=\, \prod_{i=0}^{n-1} m \big[A^{n-i} \big((A^ne^{n\varphi})^k(A^i\mathbf 1)\big)\big]^{1/n}
 \,\le\, \prod_{i=0}^{n-1} \big\| A^{n-i} \big((A^ne^{n\varphi})^k(A^i\mathbf 1)\big)\big\|^{1/n}\\[6pt]
 \,\le\, \|A\|^n \pin \big\|(A^ne^{n\varphi})^k \big\| .
\end{gather*}
Noting that
\begin{equation*}
 \big\|A_\varphi^{n(k+1)}\big\| = \sup_{m\in M(\mathcal C)}m\big[A_\varphi^{n(k+1)}\mathbf 1\big],
\end{equation*}

 \smallskip\noindent
we deduce that relations just obtained imply
\begin{equation*}
 e^{-nc}\pin\big\|A_\varphi^{n(k+1)}\big\| \le \|A\|^n\pin \big\|(A^ne^{n\varphi})^k\big\|,
\end{equation*}
and therefore
\begin{equation*}
 -nc +\ln\!\pin\big\|A_\varphi^{n(k+1)}\big\| \le n\ln\!\pin\|A\| +
 \ln\!\pin\big\|(A^ne^{n\varphi})^k\big\|.
\end{equation*}

 \medskip\noindent
Dividing the latter inequality by $k$ and turning $k\to\infty$ one gets \eqref{,,12}. \qed

\proofskip

Let us fix a measure $\mu\in M_\delta (\mathcal C)$, natural number $n$ and partition of unity $D$
in $\mathcal C$. For these objects there exists a sequence of measures $m_k\in M (\mathcal C)$ on
which the supremum in \eqref{6,,7} is attained. One may choose a subsequence $m_{k_i}$ of this
sequence such that the following limits do exist simultaneously:
\begin{equation} \label{,,13}
 \lim_{i\to\infty} m_{k_i}[A^n(g)] =: C_n(\mu,g,D), \qquad g\in D.
\end{equation}
Then by construction one has
\begin{equation} \label{,,14}
 \tau_n(\mu,D) = \sum_{g\in D}\mu[g] \ln\frac{C_n(\mu,g,D)}{\mu[g]}\pin.
\end{equation}

\begin{lemma} \label{..5}
If\/ $\tau_n(\mu,D)>-\infty$ then
\begin{equation} \label{,,15}
 \sup_{m\in M(\mathcal C)} \sum_{\substack{g\in D,\\[1pt] \mu[g]>0}} \mu[g]\,
 \frac{m[A^ng]}{C_n(\mu,g,D)} \,=\, 1.
\end{equation}
\end{lemma}

\smallskip

\proof. The finiteness of $\tau_n(\mu,D)$ and \eqref{,,14} imply that $C_n(\mu,g,D)>0$ whenever
$\mu[g]>0$. For each $m\in M (\mathcal C)$ let us consider the function
\begin{equation*}
 \eta(t) =\sum_{\substack{g\in D,\\[1pt] \mu[g]>0}} \mu[g]\ln
 \frac{(1-t)C_n(\mu,g,D) +t \, m[A^ng]}{\mu[g]}, \qquad t \in [0,1].
\end{equation*}
By definition of the numbers $C_n(\mu,g,D)$ this function attains its maximal value equal to
$\tau_n(\mu,D)$ at $t=0$. Therefore its derivative at $t=0$
\begin{equation*}
 \frac{d\eta(t)}{dt}\bigg|_{t=0} =\sum_{\substack{g\in D,\\[1pt] \mu[g]>0}}
 \mu[g]\pin \frac{m[A^ng] -C_n(\mu,g,D)}{C_n(\mu,g,D)} =
 \sum_{\substack{g\in D,\\[1pt] \mu[g]>0}} \mu[g]\pin \frac{m[A^ng]}{C_n(\mu,g,D)} -1
\end{equation*}
is nonpositive and so the left hand part in \eqref{,,15} does not exceed its right hand part.

The equality in \eqref{,,15} is attained on the sequence of measures $m_{k_i}$ from \eqref{,,13}.
\qed

\proofskip

Now we can finish the proof of Theorem~\ref{..1}.

Let us fix an arbitrary measure $\mu\in M_\delta (\mathcal C)$, natural number $n$ and partition of
unity $D$ in $\mathcal C$.

Suppose at first that there exists an element $g\in D$ satisfying the inequality $\mu[g] >0$ and
equality $m[A^n g] =0$ for all $m\in M(\mathcal C)$. Then by definition one has $\tau_n(\mu,D)
=-\infty$ and therefore $\tau(\mu) =-\infty$. Thus in this case equality \eqref{,,6} takes the form
\begin{equation} \label{,,16}
 -\infty =\inf_{\varphi\in \mathcal C} \bigl(\lambda(\varphi) -\mu[\varphi]\bigr).
\end{equation}
Let us verify it.

Consider the family of elements $\varphi_t =tg/n$, where $t\in\mathbb R$. Inequalities $0\le g\le
\mathbf 1$ and the Lagrange theorem imply that
\begin{equation*}
 e^{n\varphi_t} =e^{tg} \le \mathbf 1+e^t tg.
\end{equation*}
Therefore for each measure $m \in M(\mathcal C)$ one has
\begin{equation*}
 m[A^n(e^{n\varphi_t} \mathbf 1)] \le m[A^n(\mathbf 1+e^t tg)] =
 m[A^n\mathbf 1] + e^t t\, m[A^ng]= m[A^n\mathbf 1] \le \|A^n\|.
\end{equation*}
Thus $\| A^ne^{n\varphi_t}\| \le \|A^n\|$. Applying Lemma~\ref{..4}, we obtain the following
estimate
\begin{equation*}
 n\lambda(\varphi_t) =n\lambda(\varphi_t,A) \le \lambda(n\varphi_t,A^n) \le \ln\!\pin
 \|A^n e^{n\varphi_t}\| \le \ln\!\pin \|A^n\|.
\end{equation*}
On the other hand,
\begin{equation*}
 \mu[\varphi_t] =\mu[tg/n] = t\mu[g]/n \to +\infty \quad\text{as}\ \ t\to +\infty.
\end{equation*}
And therefore $\lambda(\varphi_t) -\mu[\varphi_t] \to -\infty$ when $t\to +\infty$. So equality
\eqref{,,16} is verified.

\medskip

It remains to consider the situation when for each element $g\in D$ satisfying the condition
$\mu[g] >0$ there exists a measure $m_g\in M(\mathcal C)$ such that $m_g[A^ng] >0$. Taking the
measure $m :=|D|^{-1}\sum_g m_g$ one obtains that
\begin{equation*}
m[A^ng] >0 \quad\text{as soon as}\ \ \mu[g]>0.
\end{equation*}
Therefore $\tau_n(\mu,D) >-\infty$. Note also that finiteness of $\tau_n(\mu,D)$ along with
\eqref{,,14} implies that the condition $\mu[g]>0$ automatically implies the inequality
$C_n(\mu,g,D)>0$.

Now let us define the family of elements
\begin{equation}\label{,,17}
\varphi_\eps :=\, \frac{1}{n}\pin\ln\Bigg\{ \sum_{\mu[g]>0} \frac{\mu[g]}{C_n(\mu,g,D)}\pin g +
\sum_{\mu[g]=0} \eps g\Bigg\}, \qquad \eps>0.
\end{equation}
For any $m\in M(\mathcal C)$ one has
\begin{align*}
 m[A^n(e^{n\varphi_\eps} \mathbf 1)] \,&=\, m\Bigg[A^n\Bigg(\sum_{\mu[g]>0}
 \frac{\mu[g]}{C_n(\mu,g,D)}\pin g\Bigg)\Bigg] +m\Bigg[A^n\Bigg(\sum_{\mu[g]=0}\eps g\Bigg)\Bigg]
 \\[6pt]
 &=\,\sum_{\mu[g]>0} \mu[g]\pin \frac{m[A^ng]}{C_n(\mu,g,D)} +
 \eps m\Bigg[A^n\Bigg(\sum_{\mu[g]=0} g \Bigg)\Bigg] \,\le\, 1 +\eps\|A^n\|
\end{align*}
(where in the final inequality we exploited Lemma~\ref{..5}).

Therefore,
\begin{equation*}
 \|A^n e^{n\varphi_\eps}\| \le 1 +\eps\|A^n\|.
\end{equation*}

 \medskip\noindent
This along with Lemma~\ref{..4} implies the estimate
\begin{equation} \label{,,18}
 n\lambda(\varphi_\eps) \le \lambda(n\varphi_\eps,A^n) \le \ln\!\pin \|A^n e^{n\varphi_\eps}\|
 \le \ln(1+\eps\|A^n\|) \le \eps\|A^n\|.
\end{equation}

On the other hand, applying concavity of logarithm and \eqref{,,14} one obtains
\begin{equation} \label{,,19}
 \begin{aligned}
   \mu[n\varphi_\eps] &\,=\,\mu\Bigg[\ln\Bigg\{\sum_{\mu[g]>0} \frac{\mu[g]}{C_n(\mu,g,D)}\pin g +
   \sum_{\mu[g]=0} \eps g\Bigg\}\Bigg] \\[6pt]
   &\,\ge\,\mu\Bigg[\sum_{\mu[g]>0} g\ln\frac{\mu[g]}{C_n(\mu,g,D)} +\sum_{\mu[g]=0} g\ln\eps\Bigg]
   =\,-\tau_n(\mu,D).
 \end{aligned}
\end{equation}
Combining \eqref{,,19} and \eqref{,,18} we get
\begin{equation*}
 \frac{\tau_n(\mu,D)}{n} \,\ge\, -\mu[\varphi_\eps] \,\ge\,
 -\mu[\varphi_\eps] +\left(\lambda(\varphi_\eps) -\frac{\eps\|A^n\|}{n}\right),
\end{equation*}
and therefore
\begin{equation*}
 \frac{\tau_n(\mu,D)}{n} +\frac{\eps\|A^n\|}{n} \,\ge\, \lambda(\varphi_\eps) -\mu[\varphi_\eps]
 \,\ge\, \inf_{\varphi\in \mathcal C}\big(\lambda(\varphi) -\mu[\varphi]\big).
\end{equation*}
This inequality along with arbitrariness of $\eps$, $n$, $D$ and definition \eqref{6,,5} of
$\tau(\mu)$ implies the inequality
\begin{equation*}
 \tau(\mu) \ge \inf_{\varphi\in \mathcal C}\big(\lambda(\varphi) -\mu[\varphi]\big).
\end{equation*}

 \medskip\noindent
Together with inequality \eqref{,,11} this proves \eqref{,,6}. \qed

\proofskip

Note now that Theorem~\ref{..2} is a straightforward corollary of the next observation.

\begin{lemma} \label{..7}
If a linear functional\/ $\mu$ on\/ $\mathcal C$ possesses the property
\begin{equation} \label{,,20}
 \inf_{\varphi\in \mathcal C} \big(\lambda(\varphi) -\mu[\varphi]\big) >-\infty,
\end{equation}
then\/ $\mu\in M_\delta (\mathcal C)$. In particular, this is true for every subgradient of the function\/
$\lambda(\varphi)$.
\end{lemma}

This lemma can be proven absolutely in the same way as the corresponding result (Lemma~7) in
\cite{BL}. The proof is based on the following properties of the functional~$\lambda(\varphi)$.

\begin{lemma} \label{..6}
The spectral potential\/ $\lambda(\varphi)$ possesses the following properties\/$:$

\smallskip

a\/$)$ if\/ $\varphi\ge\psi$, then\/ $\lambda(\varphi)\ge\lambda(\psi)$ \,\emph{(monotonicity);}

\smallskip

b\/$)$ \,$\lambda(\varphi+t) =\lambda(\varphi)+t$ for all\/ $t\in\mathbb R$ \,\emph{(additive
 homogeneity);}

\smallskip

c\/$)$ \,$|\lambda(\varphi) -\lambda(\psi)|\le \|\varphi-\psi\|$ \,\emph{(Lipschitz condition);}

\smallskip

d\/$)$ \,$\lambda((1-t)\varphi +t\psi)\le (1-t)\lambda(\varphi) +t\lambda(\psi)$ for\/ $t\in [0,1]$
\,\emph{(convexity);}

\smallskip

e\/$)$ \,$\lambda(\varphi +\delta\psi) =\lambda(\varphi+\psi)$ \,\emph{(strong
 $\delta$-invariance)}.

\end{lemma}

This lemma is proven in \cite{ABL1}, \cite{Bakhtin1}.

\begin{rem} \label{..8}
By Lemma~\ref{..6} the functional $\lambda(\varphi)$ is convex and continuous. Theorems~\ref{..1}
and \ref{..2} in essence state that the functional $-\tau(\mu)$ is the Legendre transform of
$\lambda(\varphi)$. This automatically implies that $t$-entropy $\tau(\mu)$ is concave and upper
semicontinuous (in the $^*$-weak topology) on the dual space to $\mathcal C$. In \cite{ABL1}
concavity and upper semicontinuity of $t$-entropy were proven independently and in an essentially
more complicated way.
\end{rem}

Finally we observe that variational principle for the spectral potential can be easily derived from
Theorem~\ref{..1} and Lemmas~\ref{..6}, \ref{..7}.

\begin{theorem}[\hbox spread -2pt{variational principle for the spectral potential}] \label{..10}
For each\/ $\varphi\in \mathcal C$ the following equality takes place\/$:$
\begin{equation} \label{,,21}
 \lambda(\varphi) = \max_{\mu\in M_\delta (\mathcal C)} \bigl(\tau(\mu) +\mu[\varphi]\bigr).
\end{equation}
\end{theorem}

\proof. By Lemma~\ref{..6} the functional $\lambda(\varphi)$ is convex and continuous. Thus at each
point $\varphi_0$ there exists at least one subgradient $\mu$ for $\lambda(\varphi)$. By
Lemma~\ref{..7} this subgradient belongs to $M_\delta (\mathcal C)$. By Theorem~\ref{..1} and
definition of a subgradient we have
\begin{equation*}
 \tau(\mu) =\inf_{\varphi\in \mathcal C} \big(\lambda(\varphi) -\mu[\varphi]\big) =
 \lambda(\varphi_0) -\mu[\varphi_0].
\end{equation*}
Therefore $\lambda(\varphi_0) =\tau(\mu) +\mu[\varphi_0]$. Combining this equality
with~\eqref{,,11} one obtains \eqref{,,21}. \qed

\section{Entropy Statistic Theorem}

Entropy statistic theorem is naturally formulated in terms of the dynamical system $(X, \alpha)$
corresponding to $({\mathcal C},\delta)$ (see \ref{6..4}).

Here by $M(X)$ we denote the set of all Borel probability measures on $X$. Let $x$ be an arbitrary
point of $X$. The \emph{empirical measures} $\delta_{x,n}\in M(X)$ are defined by the formula
\begin{equation} \label{9,,1}
 \delta_{x,n}(f) :=\frac{f(x)+f(\alpha(x))+\,\dotsm\,+f(\alpha^{n-1}(x))}{n} =\frac{1}{n} S_nf(x),
 \qquad f\in C(X).
\end{equation}
Evidently, the measure $\delta_{x,n}$ is supported on the trajectory of the point $x$ of length
$n$.

We endow the set $M(X)$ with the $\,^*$-weak topology of the dual space to $C(X)$. Given a measure
$\mu\in M(X)$ and its certain neighborhood $O(\mu)$ we define the sequence of sets $X_n(O(\mu))$ as
follows:
\begin{equation}\label{9,,2}
 X_n(O(\mu)) := \{\pin x\in X\mid \delta_{x,n}\in O(\mu)\pin\}.
\end{equation}

\begin{theorem}[entropy statistic theorem] \label{9..1}
Let\/ $A\!: C(X) \to C(X)$ be a certain transfer operator for\/ $(X,\alpha)$. Then for any
measure\/ $\mu\in M(X)$ and any number\/ $\eps>0$ there exist a neighborhood\/ $O(\mu)$ in the\/
$^*$-weak topology, a\/ $($large enough\/$)$ number\/~$C(\eps,\mu)$ and a sequence of functions\/
$\chi_n\in C(X)$ majorizing the index functions of the sets\/~$X_n(O(\mu))$ such that for all\/ $n$
the following estimate holds
\begin{equation}\label{,,25}
\norm{A^n\chi_n} \le C(\eps,\mu)\pin e^{n(\tau(\mu)+\eps)} .
\end{equation}
\end{theorem}

If $\tau(\mu) =-\infty$ then the number $\tau(\mu)+\eps$ in \eqref{,,25} should be replaced by
$-1/\eps$.

\medbreak

\proof. By the variational principle for $t$-entropy there exists $\varphi\in C(X)$ such that
\begin{equation*}
 \lambda(\varphi) -\mu[\varphi] < \tau(\mu)+\eps/3
\end{equation*}
(or $\lambda(\varphi) -\mu[\varphi] <-1/\eps -\eps/3$ in the case when $\tau(\mu) =-\infty)$. Let
us set
\begin{equation*}
 O(\mu) :=\big\{ \nu \in M(X) \bigm| \lambda(\varphi) -\nu [\varphi] <\tau(\mu) +\eps/3\big\}.
\end{equation*}
Then
\begin{equation*}
 X_n(O(\mu)) =\big\{ x\in X\bigm| S_n\varphi(x) =n\delta_{x,n}[\varphi] >
 n(\lambda(\varphi) -\tau(\mu) -\eps/3)\big\}.
\end{equation*}
Let
\begin{equation*}
 Y_n:= \big\{ x\in X \bigm| S_n\varphi(x) =n\delta_{x,n}[\varphi] \le
 n(\lambda(\varphi) -\tau(\mu) -\eps/2)\big\}.
\end{equation*}

 \medskip\noindent
Then by Uhryson's Lemma there exist continuous functions $\chi_n$ such that
\begin{equation*}
 0\le \chi_n \le 1,\quad \chi_n \Big(\overline{X_n(O(\mu))} \Big) =1 \ \ \text{and}\ \
 \chi_n(Y_n)=0.
\end{equation*}
Clearly, $\chi_n$ majorizes the index function of $X_n(O(\mu))$.

Take a constant $C(\eps,\mu)$ so large that
\begin{equation*}
\|A_\varphi^n\| \le C(\eps,\mu)\pin e^{n(\lambda(\varphi) +\eps/2)}, \qquad n\in\mathbb N.
\end{equation*}
Now \eqref{,,25} follows from the calculation
\begin{align*}
 C(\eps,\mu)\pin e^{n(\lambda(\varphi)+\eps/2)} \|\mathbf 1\| \,&\ge\, \|A_\varphi^n \mathbf 1\|
 \,\ge\, \| A_\varphi^n \chi_n\| \\[6pt]
 &=\, \big\|A^n\big(e^{S_n\varphi} \chi_n\big) \big\| \,\ge\,
 e^{n(\lambda(\varphi) -\tau(\mu) -\eps/2)} \| A^n \chi_n\|. \qed
\end{align*}

To summarize the material presented we recall that in \cite{ABL1} the next chain of statements for
transfer operators has been proven: `entropy statistic theorem' $\Rightarrow$ `variational
principle for the spectral potential' $\Rightarrow$ `variational principle for $t$-entropy', where
each step is rather nontrivial. In this article we obtained the inverted chain: `variational
principle for $t$-entropy' $\Rightarrow$ `variational principle for the spectral potential'
$\Rightarrow$ `entropy statistic theorem'. Thus we established equivalence: `variational principle
for $t$-entropy' $\Leftrightarrow$ `variational principle for the spectral potential'
$\Leftrightarrow$ `entropy statistic theorem'.



\end{document}